\documentclass[12pt]{article}

\usepackage{todonotes}

\usepackage{amsmath,amssymb,amsbsy,amsfonts,amsthm,latexsym,amsopn,amstext,amsxtra,euscript,amscd,color}
\usepackage{color,xcolor}
\usepackage{latexsym}
\usepackage{cite}
\usepackage{amsfonts}
\usepackage{amsfonts,mathrsfs}
\usepackage{graphicx}
\usepackage{psfrag}
\usepackage{subfigure}
\usepackage{url}
\usepackage{stfloats}
\usepackage{amsmath}
\usepackage{algorithm}
\usepackage{algorithmic}
\usepackage{hyperref}

\hypersetup{colorlinks=true}

\newcommand{\aw}[1]{{\color{magenta} #1}}

\newtheorem{theorem}{Theorem}
\newtheorem{lemma}{Lemma}

\newcommand{\quash}[1]{}

\begin{document}

\title{On the $N$th $2$-adic complexity of binary sequences identified with algebraic $2$-adic integers}

\author{Zhixiong~Chen~$^1$ and Arne~Winterhof~$^2$\\
~\\
$^1$ Fujian Key Laboratory of Financial Information Processing,\\
 Putian University, Putian, Fujian 351100, P. R. China\\
$^2$ Johann Radon Institute for Computational and Applied Mathematics,\\
Austrian Academy of Sciences, Altenbergestr.\ 69,
A-4040 Linz, Austria
}

\date{}
 
\maketitle

\begin{abstract}
 We identify a binary sequence $\mathcal{S}=(s_n)_{n=0}^\infty$ with the $2$-adic integer
  $G_\mathcal{S}(2)=\sum\limits_{n=0}^\infty s_n2^n$. In the case that $G_\mathcal{S}(2)$ is algebraic over $\mathbb{Q}$ of degree $d\ge 2$, we prove that the $N$th $2$-adic complexity of $\mathcal{S}$ is at least $\frac{N}{d}+O(1)$, where the implied constant depends only on the minimal polynomial of $G_\mathcal{S}(2)$.
  This result is an analog of the bound of M\'erai and the second author on the linear complexity of automatic sequences, that is, sequences with algebraic $G_\mathcal{S}(X)$ over the rational function field $\mathbb{F}_2(X)$. 
  
  We further discuss the most important case $d=2$ in both settings and explain that the intersection of the set of $2$-adic algebraic sequences and the set of automatic sequences is the set of (eventually) periodic sequences. Finally, we provide some experimental results supporting the conjecture that $2$-adic algebraic sequences can have also a desirable $N$th linear complexity and automatic sequences a desirable $N$th $2$-adic complexity, respectively.
  \end{abstract}

\section{Introduction}

Automatic sequences
have gained much attention during the last decades, for example, in theoretical computer science, combinatorics and number theory, see for example the monograph \cite{AS2003}.  
More recently, many famous (aperiodic) automatic sequences, including the Thue-Morse and the Rudin-Shapiro sequence, and certain subsequences such as those along squares, have been considered in view of cryptographic applications in terms of  
pseudorandomness.  
Several measures of pseudorandomness have been studied such as well distribution measure,  correlation measure of order $k$, normality, 
linear complexity, maximum order complexity
and expansion complexity, see \cite{AHN2020,DMR2019,DPS2021,GM2020,GMN2018,HW2016,M2018,MS2017,MS1998,MNW2017,MW2018CC,MW2018JNT,MW2022,P2020,S2015,SW2019UDT,SW2019,W2023}.

Automatic sequences can be defined in an algebraic way. 
For a sequence $\mathcal{S}=(s_n)_{n=0}^\infty$ over the finite field $\mathbb{F}_q$ of $q$ elements, 
 $\mathcal{S}$ is \textit{$q$-automatic} (or briefly \textit{automatic}) if and only if $G_{\mathcal{S}}(X)=\sum\limits_{n= 0}^\infty s_nX^n$ is \textit{algebraic} over $\mathbb{F}_q[X]$ by a result of Christol \cite{C1979,CKMR1980}, that is, there exists a non-zero polynomial $h(Y)\in \mathbb{F}_q[X][Y]$ (with coefficients in $\mathbb{F}_q[X]$) such that 
$$h(G_{\mathcal{S}}(X))=0,$$ 
see also \cite[Theorem 12.2.5]{AS2003}.  Sometimes, an automatic sequence is also called an \emph{expansion sequence} \cite{D2012,GMN2018}.  

Automatic sequences guarantee some desirable features of pseudorandomness.
In particular, M\'{e}rai and the second author proved in \cite{MW2018JNT} that any aperiodic automatic sequence over $\mathbb{F}_q$ has $N$th linear complexity of order of magnitude $N$. 
In particular, they proved the lower bound $\frac{N}{\deg(h)}+O(1)$,
where the implied constant depends only on $h(Y)$ and $\deg(h)$ is the degree of $h(Y)$. 
In the best case, that is $\deg(h)=2$, the $N$th linear complexity is $\frac{N}{2}+O(1)$,
which is essentially the expected value of the $N$th linear complexity of a random sequence \cite{N1990}.

In this work, we study analogs of automatic sequences in the $2$-adic setting and prove similar bounds on their $2$-adic complexity. 
In particular, we restrict ourselves to binary sequences 
$\mathcal{S}=(s_n)_{n=0}^\infty$ with $s_n\in\{0,1\}$. 
The quantity 
$$G_{\mathcal{S}}(2)=\sum_{n=0}^\infty s_n2^n$$ 
is a $2$-adic integer and plays an important role in our study.

We denote by $\mathbb{Z}_2$ the set of {\em $2$-adic integers}   
$$
\sum\limits_{n=0}^{\infty}a_n2^n ~~~\mathrm{with}~~~a_n\in\{0,1\}.
$$
Note that these series converge with respect to the {\em $2$-adic absolute value}
$$\left|\sum_{n=k}^\infty a_n 2^n\right|_2=2^{-k}\rightarrow 0\quad \mbox{if }a_k=1,\quad k\rightarrow \infty.$$
The fraction field of $\mathbb{Z}_2$, denoted by $\mathbb{Q}_2$, consists of all infinite sums
$$
\sum\limits_{n=k}^{\infty}a_n2^n ~~~\mathrm{with}~~~a_n\in\{0,1\}, ~~k\in\mathbb{Z},
$$
with finitely many terms of negative degree and possibly infinitely many terms of
non-negative degree.   

An element $\alpha\in\mathbb{Q}_2$ is \textit{algebraic} over $\mathbb{Z}$ 
if there is a non-zero polynomial
 $h(Y)\in \mathbb{Z}[Y]$ such that $h(\alpha)=0$.
If $G_\mathcal{S}(2)$ is algebraic, we also call the sequence $\mathcal{S}$ a {\em $2$-adic algebraic sequence}.
For example, any rational number $\frac{a}{b}$ with $a,b\in\mathbb{Z}$, $b\not=0$, is algebraic since we can take $h(Y)=bY-a$.
It is obvious that, if $\mathcal{S}$ is eventually periodic, then $G_{\mathcal{S}}(2)$ is rational and thus algebraic over~$\mathbb{Z}$, see for example \cite[Theorem~4.2.4]{GK2012}.

Now we introduce the $2$-adic complexity, which is the analog of the linear complexity in the $2$-adic setting.
Goresky and Klapper proposed a new type of registers called the feedback with carry shift registers (FCSRs) in the 1990s. The FCSRs have the advantage of simple hardware implementation in high speed and potential applications in cryptography and communications. For more details about FCSRs we refer to the monograph~\cite{GK2012}.

A binary FCSR of length $L$ is carried out to produce a binary sequence $\mathcal{S}=(s_n)_{n=0}^\infty$ over $\{0,1\}$ by  a recurrence relation in the ring $\mathbb{Z}$ of integers with $a_j\in\{0,1\}$ for $0\leq j\leq L-1$,
$$
s_{i+L}+2z_{i+L}=a_{L-1}s_{i+L-1}+a_{L-2}s_{i+L-2}+\ldots+a_0s_{i}+z_{i+L-1}, \quad i\geq 0,
$$
with 
initial integer $z_{L-1}$ and carries $z_{L},z_{L+1},\ldots$ 
The notion of \textit{$N$th $2$-adic complexity}  of $\mathcal{S}$,  denoted by
$\lambda_{\mathcal{S}}(N)$,  is raised to 
approximate the length of a shortest FCSR that generates the  first $N$ elements of $\mathcal{S}$.
$\lambda_{\mathcal{S}}(N)$ is the binary logarithm of
$\Lambda_{\mathcal{S}}(N)$ with 
$$
\Lambda_{\mathcal{S}}(N)=\min \left\{\max\{|f|, q\} : f,q\in\mathbb{Z}, q>0 \mbox{ odd},~ q\sum_{i=0}^{N-1}s_i2^i \equiv f \bmod 2^N   \right\};
$$
see \cite[p.328]{GK2012} or the recent survey \cite{W2023}, that is, $\Lambda_{\mathcal{S}}(N)=2^{\lambda_{\mathcal{S}}(N)}$. 
Here we use the notation
$$a\equiv b\bmod c\quad \mbox{if and only if}\quad a=b+kc \mbox{ for some }k\in \mathbb{Z}$$
for any integers $a,b,c$.
Any pair $(q,f)$ 
with $\Lambda_{\mathcal{S}}(N)=\max\{|f|,q\}$ is called a \textit{minimal rational representation} for the first $N$ elements of $\mathcal{S}$.
It is trivial that $\lambda_{\mathcal{S}}(N)\leq \lambda_{\mathcal{S}}(N+1)$ for $N\geq 1$.

In Section \ref{Sect-2adic} we prove bounds on the $N$th $2$-adic complexity of a binary sequence~$\mathcal{S}$ if $G_\mathcal{S}(2)\in \mathbb{Z}_2$ is algebraic over $\mathbb{Z}$. 
In Section \ref{sect-roots}, with the help of Hensel’s lemma, we discuss the existence of zeros in $\mathbb{Z}_2$ of polynomials of degree $2$ defined over $\mathbb{Z}$ and of zeros in $\mathbb{F}_2[[X]]$ of polynomials of degree $2$ over $\mathbb{F}_2[X]$, respectively.
 In the final section, we discuss that the sets of non-periodic automatic and non-periodic $2$-adic algebraic sequences are disjoint as well as we discuss whether automatic sequences can have a desirable $N$th $2$-adic complexity and $2$-adic algebraic sequences can have desirable $N$th linear complexity as well. Moreover, we mention a consequence of Roth's theorem to the $2$-adic complexity.

\section{$N$th $2$-Adic complexity of aperiodic sequences}\label{Sect-2adic} 

In this section we prove the following analog of \cite[Theorem~1]{MW2018JNT} on the linear complexity of automatic sequences for the $2$-adic complexity of sequences with algebraic~$G_\mathcal{S}(2) \in \mathbb{Z}_2$.

\begin{theorem}\label{Thm-bound}
Let $\mathcal{S}=(s_n)_{n=0}^\infty$ be  
an aperiodic, (that is, not eventually periodic) sequence  over~$\{0,1\}$ and $$h(Y)=h_0+h_1Y+\ldots+h_dY^d\in \mathbb{Z}[Y]$$
be a polynomial of degree $d$ without rational roots such that $h(G_{\mathcal{S}}(2))=0$. Put
$$H=\sum_{i=0}^d|h_i|.$$ 
Then for $N\geq 1$ we have
$$
\frac{N}{d}-\frac{\log_2(H)}{d}\leq \lambda_{\mathcal{S}}(N)\leq \frac{d-1}{d}N+\frac{\log_2(H)}{d}+1.3.
$$
\end{theorem}

Proof. Since $\mathcal{S}$ is not eventually periodic, $G_{\mathcal{S}}(2)$ is not a rational number and thus~$d\geq 2$. 

Assume that $(q,f)$ is a minimal rational representation for the first $N$ elements of $\mathcal{S}$, that is,
$$
qG_{\mathcal{S}}(2)\equiv f \bmod 2^N
$$
and
$$\Lambda_\mathcal{S}(N)=\max\{q,|f|\}.$$
Then we have
$$
0=q^d\cdot h(G_{\mathcal{S}}(2))\equiv \sum\limits_{i=0}^{d}h_if^iq^{d-i}  \bmod{2^N}.
$$
We remark that $\sum\limits_{i=0}^{d}h_if^iq^{d-i}\neq 0$ since otherwise $\frac{f}{q}$ is a rational root of $h(x)$, which contradicts our assumption.  
Hence
$$
2^N\leq \left|\sum\limits_{i=0}^{d}h_if^iq^{d-i}\right|\leq H\cdot \Lambda_{\mathcal{S}}(N)^{d},
$$
which implies the lower bound on $\lambda_{\mathcal{S}}(N)$.   

We prove the upper bound  by induction on $N$. It is trivial for $N=1$ since $\lambda_{\mathcal{S}}(1)=0$. 

Assume that the upper bound holds for $N$ and we estimate $\lambda_{\mathcal{S}}(N+1)$. With the lower bound on $\lambda_{\mathcal{S}}(N)$ and \cite[Lemma 18.5.1]{GK2012} we have with $\varepsilon=1.3$,
\begin{eqnarray*}
\lambda_{\mathcal{S}}(N+1) &\leq & \log_2\left(2^{\lambda_{\mathcal{S}}(N)}+2^{N-\lambda_{\mathcal{S}}(N)}\right)\\
&\leq & \log_2\left(2^{\frac{d-1}{d}N+\frac{\log_2(H)}{d}+\varepsilon}+2^{\frac{d-1}{d}N+\frac{\log_2(H)}{d}}\right)\\
&\leq & \frac{d-1}{d}(N+1)+\frac{\log_2(H)}{d}+\varepsilon+\log_2\left(2^{-\frac{d-1}{d}}+2^{-\frac{d-1}{d}-\varepsilon}\right)\\
&< & \frac{d-1}{d}(N+1)+\frac{\log_2(H)}{d}+\varepsilon,     
\end{eqnarray*}
where we used $\log_2\left(2^{-\frac{d-1}{d}}+2^{-\frac{d-1}{d}-\varepsilon}\right)<\log_2(1)=0$, 
which finishes the proof.
\qed\\

%We note that $\varepsilon>v_d:=\log_2\left(\frac{1}{2^{\frac{d-1}{d}}-1}\right)$, where $v_2\approx 1.271$, $v_3\approx 0.767$, $v_4\approx 0.552$, $v_5\approx 0.432$, $v_6\approx 0.355$. 
 
In \cite{CW2024} the authors proved that, for fixed $N$, the expected value of the $N$th $2$-adic complexity over all binary sequences of length $N$ is  
close to $\frac{N}{2}$ and the deviation from~$\frac{N}{2}$ is at most of order of magnitude $\log(N)$. 
Moreover, for a random binary sequence $\mathcal{S}$, and running $N$, the $N$th $2$-adic complexity satisfies with probability~1
$$
\lambda_{\mathcal{S}}(N)=\frac{N}{2}+O(\log(N)) \quad \mbox{for all }N=1,2,\ldots
$$ 
Hence, with respect to the $2$-adic complexity, any $\mathcal{S}$ for which $G_{\mathcal{S}}(2)$ is a root of an irreducible polynomial $h(Y)\in\mathbb{Z}[Y]$ of degree $d=2$ behaves like a random sequence.

For the simplest case $h(Y)=Y^2-D\in\mathbb{Z}[Y]$,  
from \cite[Thm.2, Sect.\ 1.6, p.49]{BS1986} 
we know that $h(Y)$ has a zero in $\mathbb{Z}_2$ if and only if $D$ is of the form 
$$D=4^{m}d, \quad m\geq 0,\quad d\equiv 1 \bmod{8},$$
or $D=0$.
By Theorem~\ref{Thm-bound} we have 
$$\frac{N-\log_2(|D|+1)}{2}\le \lambda_\mathcal{S}(N)\le \frac{N+\log_2(|D|+1)}{2}+1.3.$$

For example,  the binary expansion of $\sqrt{17}$ is 
$$
\sqrt{17}=1 + 2^3 + 2^5 + 2^6 + 2^7 + 2^9 + O(2^{10})\in\mathbb{Z}_2,
$$
where $O(2^{10})$ denotes a term of the form $A=\sum\limits_{n=10}^\infty a_n2^n$, that is,
$|A|_2\le 2^{-10}$.
The binary expansion 
of $-\sqrt{17}$ is exactly the complement of the binary expansion of $\sqrt{17}$ except the first term since 
$$\sum_{n=0}^\infty 2^n=-1.$$
For $D=-7$, we get
$$
\sqrt{-7}=1 + 2^2 + 2^4 + 2^5 + 2^7 + O(2^{10})\in\mathbb{Z}_2.
$$\\

We discuss the general case of polynomials of degree $2$ in the following section.

\section{Roots of polynomials of degree $2$}\label{sect-roots}

In this section we discuss the zeros of quadratic equations in more detail, which provide the strongest bounds $\frac{N}{2}+O(1)$ on the $N$th linear complexity for non-rational zeros $G_\mathcal{S}(X)\in \mathbb{F}_2[[X]]$ over $\mathbb{F}_2(X)$ and on the $N$th $2$-adic complexity for non-rational zeros $G_\mathcal{S}(2)\in \mathbb{Z}_2$ over $\mathbb{Q}$, respectively.

\subsection{Roots in $\mathbb{Z}_2$}%\label{Sect-roots-inZ2}

Let 
$$h(Y)=aY^2+bY+c\in\mathbb{Z}[Y],\quad a\not=0,$$ 
be a polynomial of degree $2$ with integer coefficients. Now we study the zeros of $h(Y)$ in $\mathbb{Z}_2$.
Since $h(Y)$ and $d\cdot h(Y)$ have the same zeros in $\mathbb{Q}_2$ for any integer $d\neq 0$, we may restrict ourselves to the case 
\begin{equation}\label{gcd}\gcd(a,b,c)=1.
\end{equation}
Such $h(Y)$ are called \textit{primitive polynomials} over $\mathbb{Z}$, see for example \cite[Sect.11.2, p.181]{I2004}.

We assume that $h(Y)$ is irreducible over~$\mathbb{Q}$ since otherwise it has two rational zeros
and the corresponding sequences are eventually periodic.
We have to distinguish the parities of the coefficients of $h(Y)$.\\

\underline{Case 1}: $b$ is even.

\underline{Case 1.1}: $a$ is odd. \\
In this case $a$ is invertible in $\mathbb{Z}_2$. Hence, $h(Y)$ has two zeros in $\mathbb{Z}_2$ if and only if the  discriminant $b^2-4ac=4^m D$ with $D\equiv 1\bmod{8}$, that is, $\sqrt{D}\in \mathbb{Z}_2$. In this case $m\ge 1$ and the zeros are 
$$ \frac{-\frac{b}{2}\pm 2^{m-1}\sqrt{D}}{a}. $$
For example, take $h(Y)=3Y^2-4Y+9$, then the  two zeros are
$$
\frac{2 + \sqrt{-23}}{3}=1 + 2^2 + 2^5 + 2^6 + 2^8 + O(2^{10})\in\mathbb{Z}_2
$$
and
$$\frac{2 - \sqrt{-23}}{3}=1 + 2 + 2^2 + 2^6 + 2^8+O(2^{10})\in\mathbb{Z}_2.$$

\underline{Case 1.2}: $a$ is even.\\
The case $c\equiv 0 \bmod 2$ is excluded by \eqref{gcd}.\\ %since we assume $\gcd(a,b,c)=1$.\\
If $c\equiv 1 \bmod 2$, we have $h(Y)\equiv 1 \not\equiv 0\bmod{2}$, that is, $h(Y)$ has no zero modulo $2$ and thus no zero in $\mathbb{Z}_2$.\\

For the case of odd $b$, Hensel's lemma 
is a strong tool to prove the existence of zeros in $\mathbb{Z}_2$ of rational polynomial equations. We state it here for the convenience of the reader, see for example \cite[Theorem 3, p.16]{K1984}. 

\begin{lemma}\label{Hensel}
For a prime $p$ let $\mathbb{Z}_p$ be the ring of $p$-adic integers. Write
$$F(Y)=c_0+c_1Y+\ldots+c_nY^n\in \mathbb{Z}_p[Y]$$ and denote by $F'(Y)$
the derivative of $F(Y)$. If 
$$F(\alpha_0)\equiv 0\bmod{p}\quad \mbox{and}\quad F'(\alpha_0)\not\equiv 0\bmod{p}$$ 
for some $\alpha_0\in\mathbb{Z}_p$, then there exists a unique $p$-adic integer $\alpha\in\mathbb{Z}_p$ such that 
$$F(\alpha)=0 ~\mbox{ and } \alpha\equiv \alpha_0 \bmod{p}. $$
Here for two $p$-adic integers $\alpha=\sum\limits_{i=0}^\infty a_ip^i$ and $\beta=\sum\limits_{i=0}^\infty b_ip^i$ we used the notation 
$$\alpha\equiv \beta\bmod p\quad \mbox{if and only if}\quad a_0=b_0.$$
\end{lemma}

Now we continue the discussion.

\underline{Case 2}: $b$ is odd.\\
We see that the derivative of $h(Y)$ satisfies
\begin{equation}\label{derivnot0}
h'(Y)=2aY+b\equiv 1\bmod 2.  
\end{equation}

\underline{Case 2.1}: $a$ is odd.

\underline{Case 2.1.1}: $c$ is odd.\\ 
Since $h(Y)\equiv Y^2+Y+1 \bmod 2$ and thus $h(0)\equiv h(1)\equiv 1 \bmod 2$, $h(Y)$ has no zero modulo $2$ and thus no zero in~$\mathbb{Z}_2$.

\underline{Case 2.1.2}: $c$ is even.\\
Since $h(Y)\equiv Y^2+Y \bmod 2$, we have $h(0)\equiv h(1)\equiv 0 \bmod 2$ and thus $h(Y)$ has two different zeros in $\mathbb{Z}_2$ by~\eqref{derivnot0} and Lemma \ref{Hensel}.\\
For example, take $h(Y)=Y^2+5Y+2$, then the two zeros in $\mathbb{Z}_2$ are
$$\frac{-5+\sqrt{17}}{2}=2^1+2^4+2^5+2^6+2^8+2^9+O(2^{10})\in\mathbb{Z}_2$$
and $$\frac{-5-\sqrt{17}}{2}=1+2^3+2^7+O(2^{10})\in\mathbb{Z}_2.$$

\underline{Case 2.2}: $a$ is even.\\
We have $h(Y)\equiv Y+c \bmod{2}$, which means that 
$$h(c+1)\not\equiv h(c)\equiv 0 \bmod{2}$$ 
and thus $h(Y)$ has one zero in $\mathbb{Z}_2$ by \eqref{derivnot0} and Lemma \ref{Hensel}.\\
For example, consider $h(Y)=2Y^2+Y+1$. Then we see that there is one zero in $\mathbb{Z}_2$ but the other zero is in $\mathbb{Q}_2\setminus \mathbb{Z}_2$:
$$ \frac{-1+\sqrt{-7}}{4}=1+2^2+2^3+2^5+ O(2^{10})\in\mathbb{Z}_2$$
 but 
 $$\frac{-1-\sqrt{-7}}{4}=2^{-1} + 2 + 2^4 + 2^6 + 2^7 + 2^8 + 2^9 + O(2^{10})
 \not\in\mathbb{Z}_2.$$

\subsection{Roots in $\mathbb{F}_2[[X]]$}

We recall that in \cite[Theorem 1]{MW2018JNT} for aperiodic $q$-automatic sequences, a lower and an upper bound on the $N$th linear complexity have been proved via the polynomial $h(Y)\in \mathbb{F}_q[X][Y]$ with no rational roots and with $h(G_{\mathcal{S}}(X))=0$. As mentioned before, if $\deg(h)=2$, the $N$th linear complexity is  close to $\frac{N}{2}$ and $\mathcal{S}$ behaves like a random sequence in terms of the $N$th linear complexity. So it is interesting to discuss the zeros in $\mathbb{F}_2[[X]]$ of the quadratic polynomial 
$$h(Y)=a(X)Y^2+b(X)Y+c(X)$$ 
with $a(X),b(X),c(X) \in\mathbb{F}_2[X]$ and 
\begin{equation}\label{gcdpol}
\gcd(a(X),b(X),c(X))=1.
\end{equation}
It is clear that
 \begin{equation}\label{eq-hY}
     h(Y)\equiv a(0)Y^2+b(0)Y+c(0)\bmod {X},
 \end{equation}
from which we see that 
$h(0)\equiv c(0) \bmod{X}$ and $h(1)\equiv a(0)+b(0)+c(0) \bmod{X}$, so we distinguish the triples $(a(0),b(0),c(0))\in \mathbb{F}_2^3$. 

Similarly as above we will use  Hensel's lemma for formal power series, see for example \cite[Lemma 4.6, p.129]{N1999}.

\begin{lemma}\label{Hensel-poly} 
Let $\mathbb{F}_2[[X]]$ be the local ring of formal power series over $\mathbb{F}_2$. 
Write
$$F(Y)=c_0(X)+c_1(X)Y+\ldots+c_n(X)Y^n\in \mathbb{F}_2[[X]][Y]$$ and let $F'(Y)=c_1(X)+2c_2(X)Y+\ldots+nc_n(X)Y^{n-1}$ 
be the formal derivative of $F(Y)$. 

If 
$$F(\alpha_0(X))\equiv 0 \bmod X \quad \mbox{and}\quad F'(\alpha_0(X))\not\equiv 0\bmod X$$ 
for some $\alpha_0(X)\in\mathbb{F}_2[[X]]$, then there exists a unique $\alpha(X)\in\mathbb{F}_2[[X]]$ such that 
$$F(\alpha(X))=0 \quad \mbox{and} \quad \alpha(X)\equiv \alpha_0(X) \bmod X. $$
Note that we use 
$$\sum_{i=0}^\infty a_iX^i\equiv \sum_{i=0}^\infty b_iX^i\bmod X^n
\quad\mbox{if and only if}\quad (a_0,\ldots,a_{n-1})=(b_0,\ldots,b_{n-1}).$$
\end{lemma}

Now we turn to consider the zeros of $h(Y)$.  

\underline{Case 1}: $(a(0),b(0),c(0))=(0,0,0)$.\\
This case is excluded by \eqref{gcdpol}.

\underline{Case 2}: $(a(0),b(0),c(0))\in \{(0,0,1),(1,1,1)\}$.\\
$h(Y)$ has no zero in $\mathbb{F}_2[[X]]$ by~\eqref{eq-hY} since $h(0)\equiv h(1)\equiv 1 \bmod{X}$.

\underline{Case 3}: $(a(0),b(0),c(0))\in \{(0,1,0),(0,1,1)\}$.\\
 $h(Y)$ has one zero in $\mathbb{F}_2[[X]]$ 
 since $h(c(0)+1)\not\equiv h(c(0))\equiv 0 \bmod X$ and  $h'(Y)\equiv b(0)\equiv 1 \bmod{X}$
  by Lemma \ref{Hensel-poly}.\\
For example,  the zeros of $h(Y)=XY^2+Y+X$ are
$$
Y_1=X+X^3+X^7+O(X^{10}) \in \mathbb{F}_2[[X]]
$$
and
$$
Y_2=X^{-1}+X+X^3+X^7+O(X^{10}) \not\in \mathbb{F}_2[[X]],
$$
where $O(X^{10})$ denotes a term of the form $A=\sum\limits_{n=10}^\infty a_nX^n$
of absolute value $|A|_X\le 2^{-10}$.

The zeros of $h(Y)=XY^2+Y+1$ are
$$
Y_1=1+X+X^3+X^7+O(X^{10}) \in \mathbb{F}_2[[X]]
$$
and
$$
Y_2=X^{-1}+1+X+X^3+X^7+O(X^{10}) \not\in \mathbb{F}_2[[X]].
$$

\underline{Case 4}: $(a(0),b(0),c(0))=(1,1,0)$.\\  
 $h(Y)$ has two zeros in $\mathbb{F}_2[[X]]$ by \eqref{eq-hY} and Lemma \ref{Hensel-poly} with $h'(Y)\equiv b(0)\equiv 1 \bmod{X}$, since $h(0)\equiv h(1)\equiv 0 \bmod{X}$.\\
For example,  the zeros of $h(Y)=(1+X)Y^2+Y+X^2$ are
$$
Y_1=X^2+X^4+X^5+X^8+X^9+O(X^{10}) \in \mathbb{F}_2[[X]]
$$
and
$$
Y_2=1+X+X^3+X^6+X^7+O(X^{10}) \in \mathbb{F}_2[[X]].
$$

\underline{Case 5}: $(a(0),b(0),c(0))\in \{(1,0,1),(1,0,0)\}$.

\underline{Case 5.1}: $b(X)\not=0$.\\ 
Write 
$$b(X)=X^m b_1(X),\quad m\ge 1,\quad \gcd(b_1(X),X)=1.$$
Then modulo $X^m$ the equation $h(Y)=0$ simplifies to 
$$h(Y)\equiv a(X)Y^2+c(X)\equiv 0\bmod X^m.$$
A solution modulo $X^m$ exists if and only if $c(X)a(X)^{-1}\bmod X^m$ is a square, that is,
$$c(X)a(X)^{-1}\equiv d_1(X)^2\equiv d_1(X^2)\bmod X^m$$
for some $d_1(X)\in \mathbb{F}_2[X]$.
Now we make a change of variables $Y$ to $Z$ by substituting
$$Y=X^mZ+d_1(X)$$
and get with
$$d(X)=\frac{a(X)d_1(X)^2+c(X)}{X^m}+b_1(X)d_1(X)\in \mathbb{F}_2[X]$$
after canceling $X^m$
$$\overline{h}(Z)=X^{-m}h(X^mZ+d_1(X))=a(X)X^{m}Z^2+X^{m}b_1(X)Z+d(X).$$
Note that there is a one-to-one correspondence between the roots of $h(Y)$ and $\overline{h}(Z)$ but $\overline{h}(Z)$ belongs to a different case.

If $d(X)\not\equiv 0\bmod X^m$, then dividing by $\gcd(d(X),X^m)$ we reduce this case to
\begin{itemize}
\item Case 2 ($(a(0),b(0),c(0))=(0,0,1)$) with no zero
\end{itemize}
and otherwise by dividing by $X^m$ either to 
\begin{itemize}
\item Case 4 ($(a(0),b(0),c(0))=(1,1,0)$) with two zeros or to 
\item Case 2 ($(a(0),b(0),c(0))=(1,1,1)$) with no zero.
\end{itemize}
 For example, we take $h(Y)=(1+X^2)Y^2+X^5Y+X^2$ and reduce to $h(Y)\equiv (1+X^2)Y^2+X^2 \bmod{X^5}$. Then we get
 $$
 d_1(X)^2=d_1(X^2)\equiv \left(\frac{X}{1+X}\right)^2\equiv X^2+X^4 \bmod X^5,
 $$
$d(X)=X^2$ and
  $$
\overline{h}(Z)=X^{-5}h(X^5Z+X+X^2)=(1+X^2)X^5Z^2+X^5Z+X^2.
 $$
 Dividing by $X^2$ we see that $X^{-2}\overline{h}(Z)$ is of the form of
Case 2 and thus $h(Y)$ has no zero.\\

\underline{Case 5.2}: $b(X)=0$.\\
If $a(X)=A(X)^2=A(X^2)$ and $c(X)=C(X)^2=C(X^2)$ for some $A(X),C(X)\in\mathbb{F}_2[X]$, then $h(Y)$ has a (unique) zero $C(X)/A(X)$, which is rational and corresponds to an eventually periodic sequence. If no such $A(X)$ or $C(X)$ exists, then $h(Y)$ has no zero. \\ 
For example, take $h(Y)=(1+X^2+X^4)Y^2+X^6$. Then the zero is 
$$
\frac{X^3}{1+X+X^2}=\frac{X^3(1+X)}{1+X^3}=(X^3+X^4)\sum\limits_{n=0}^{\infty}X^{3n}
 \in \mathbb{F}_2[[X]].
$$\\

Finally we take the Thue-Morse sequence as an example of an algebraic sequence of degree $2$. 
The Thue-Morse sequence $\mathcal{T}=(t_n)_{n=0}^\infty$ over $\mathbb{F}_2=\{0,1\}$ is defined by
$$
t_n=
\left\{
\begin{array}{cc}
t_{n/2},       & n ~~\mbox{even},\\
1-t_{(n-1)/2}, & n ~~\mbox{odd},
\end{array}
\right.
\quad n=1,2,\ldots
$$
with initial value $t_0=0$, 
 see for example \cite{AS2003,CKMR1980,MW2018CC}.  
The generating function  $G_{\mathcal{T}}(X)$ satisfies $h(G_{\mathcal{T}}(X))=0$ with
$$
h(Y)=(X+1)^3Y^2 + (X+1)^2Y+X \in\mathbb{F}_2[X][Y],
$$
which has two zeros in $\mathbb{F}_2[[X]]$ (Case 4):
$$G_\mathcal{T}(X)=X+X^2+X^4+X^7+X^8+O(X^{11})$$ 
and 
$$G_{\mathcal{T}'}(X)=G_\mathcal{T}(X)+\frac{1}{X+1}=G_\mathcal{T}(X)+\sum_{n=0}^\infty X^n=1+X^3+X^5+X^6+X^9+X^{10}+O(X^{12}),$$
which can be identified with the dual sequence $\mathcal{T}'=(1-t_n)_{n=0}^\infty$ of the Thue-Morse sequence.
By \cite[Theorem 1]{MW2018JNT}, the $N$th linear complexities $L_{\mathcal{T}}(N)$ \mbox{and $L_{\mathcal{T}'}(N)$ are} bounded by
\begin{equation}\label{TM}
\left\lceil\frac{N-1}{2}\right\rceil\leq L_{\mathcal{T}}(N),~ L_{\mathcal{T}'}(N)\leq \left\lfloor\frac{N+2}{2}\right\rfloor.
\end{equation}
More precisely, M\'erai and the second author also showed %in \cite{MW2018JNT}
that the exact value  of the~$N$th linear complexity of the Thue Morse sequence is 
$$
L_{\mathcal{T}}(N)=2\left\lfloor\frac{N+2}{4}\right\rfloor
 ~~ \mbox{ for all } N\ge 1$$
 by considering  the continued fraction of $G_{\mathcal{T}}(X^{-1})$, see \cite[Theorem 2]{MW2018JNT} or the recent survey \cite[Theorem 3.4]{MW2022}. For more details see also the master thesis \cite{G2024}, where the dual of the Thue-Morse sequence is studied, as well,
$$L_{\mathcal{T}'}(N)=2\left\lfloor \frac{N}{4}\right\rfloor+1  ~~ \mbox{ for all } N\ge 1.$$
Hence, for any $N$ one of the sequences $\mathcal{T}$ or $\mathcal{T}'$ attains the lower bound of \eqref{TM} and the other one the upper bound.

\section{Final remarks}

\subsection{Automatic sequences vs.\ $2$-adic algebraic sequences}
In this work, we studied analogs of automatic sequences in the $2$-adic setting to prove a lower and an upper bound on the $N$th $2$-adic complexity of an aperiodic binary sequence $\mathcal{S}$ if the $2$-adic integer $G_{\mathcal{S}}(2)$ is algebraic over $\mathbb{Q}$.

A similar result for the linear complexity of non-periodic automatic, that is irrational $\mathbb{F}_2(X)$-algebraic, sequences was proved in \cite{MW2018JNT}.
This may raise the natural question whether there are sequences which are both automatic and $2$-adic algebraic.

It is mentioned in \cite[p.\ 930/931]{BY2017} that following the proof of \cite[Theorem~6]{AB2007} one can show that the $2$-adic number $G_\mathcal{S}(2)$ is transcendental if $\mathcal{S}$ is a non-periodic $2$-automatic sequence,  see also \cite[Theorems 2 and 6]{AB2007}, \cite[Theorem AB]{ACG2020} or \cite{BY2017} for more detailed discussions.
 Hence, these two sets of binary sequences are disjoint. Moreover, a sequence $\mathcal{S}$ is  eventually periodic if and only if $G_{\mathcal{S}}(X)$ is algebraic over $\mathbb{F}_2[X]$ and  $G_\mathcal{S}(2)$ is algebraic over $\mathbb{Q}$,  see \cite[Theorem 7]{AB2007}. In this case, $G_\mathcal{S}(X)$ is a rational function in $\mathbb{F}_2(X)$ as well as $G_\mathcal{S}(2)$ is a rational number in $\mathbb{Q}$.

 \subsection{Linear complexity vs.\ $2$-adic complexity}
Consequently, we can neither apply Theorem \ref{Thm-bound} to bound the $N$th $2$-adic complexity of aperiodic $2$-automatic sequences nor \cite[Theorem~1]{MW2018JNT} to estimate the $N$th linear complexity of $2$-adic algebraic sequences.  
However, there is no reason why an automatic sequence cannot have also a desirable $N$th $2$-adic complexity or why a $2$-adic algebraic sequence cannot have also a desirable $N$th linear complexity.

For example, for the Thue-Morse sequence $\mathcal{T}$, it was proved in a different way that 
$$\lambda_{\mathcal{T}}(N)\ge \frac{N}{5} \quad \mbox{for }N\ge 4$$ 
in \cite[p.6062]{CCOW2024} and 
$$\frac{N}{4}-3<\lambda_{\mathcal{T}}(N)< \frac{3N}{4}+1\quad \mbox{for }N\ge 1$$ 
in the master thesis \cite[Theorem 5.1.3]{O2023}. Experimental results in \cite[Sect.\ 5]{O2023} indicate that $\lambda_{\mathcal{T}}(N)$ is very close to $\frac{N}{2}$, which supports the conjecture that the Thue-Morse sequence behaves like a truly random sequence 
not only in terms of the $N$th linear complexity but also in terms of the $N$th $2$-adic complexity, see \cite{CW2024}. 

\begin{figure}[hp]
\begin{center}
\includegraphics[scale=0.32]{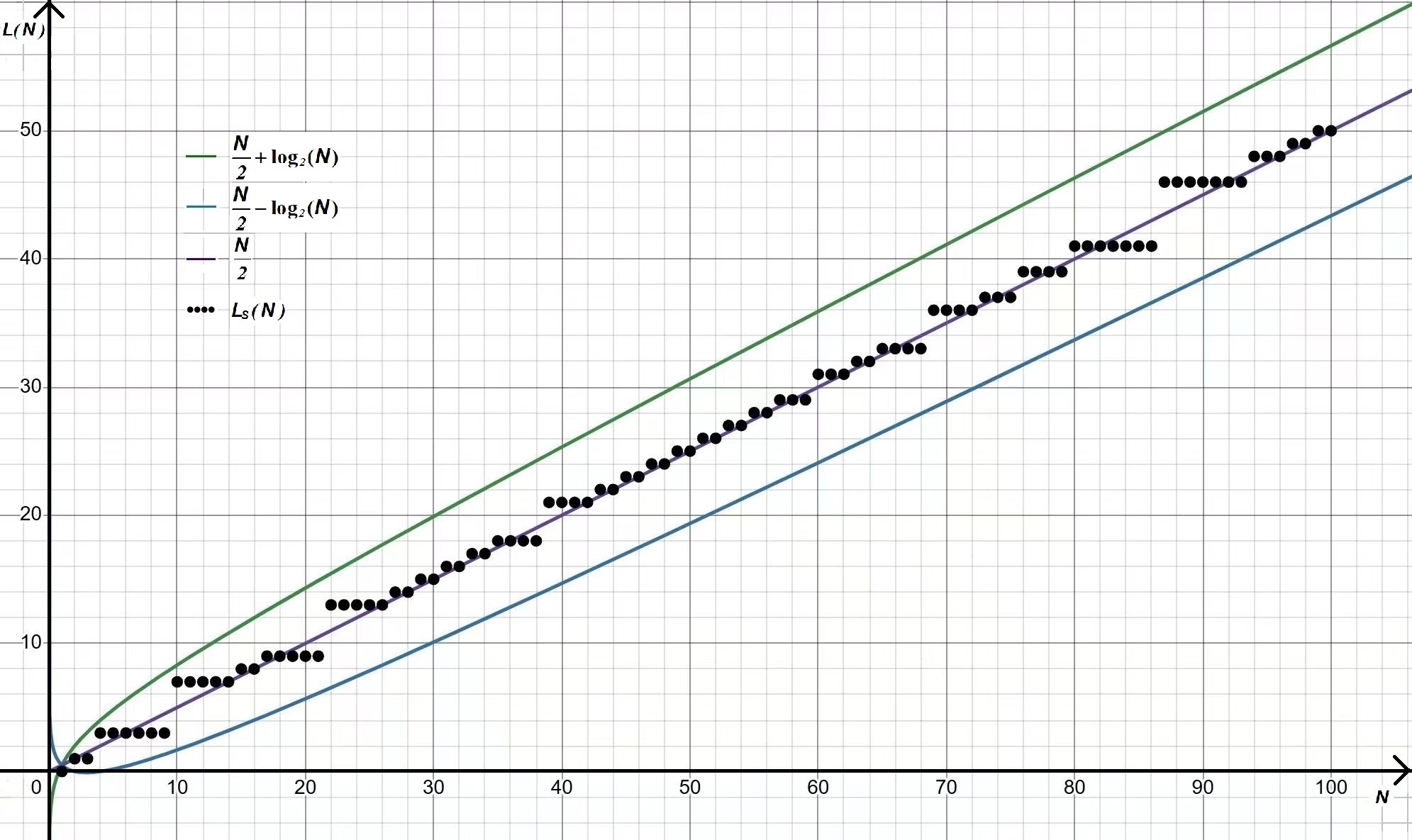} 
\end{center}
\caption{The $N$th linear complexity of $\mathcal{S}$ with $G_{\mathcal{S}}(2)=\sqrt{17}$ for $N\le 100$}
\label{fig=17}
\end{figure}

\begin{figure}[!h]
\begin{center}
\includegraphics[scale=0.32]{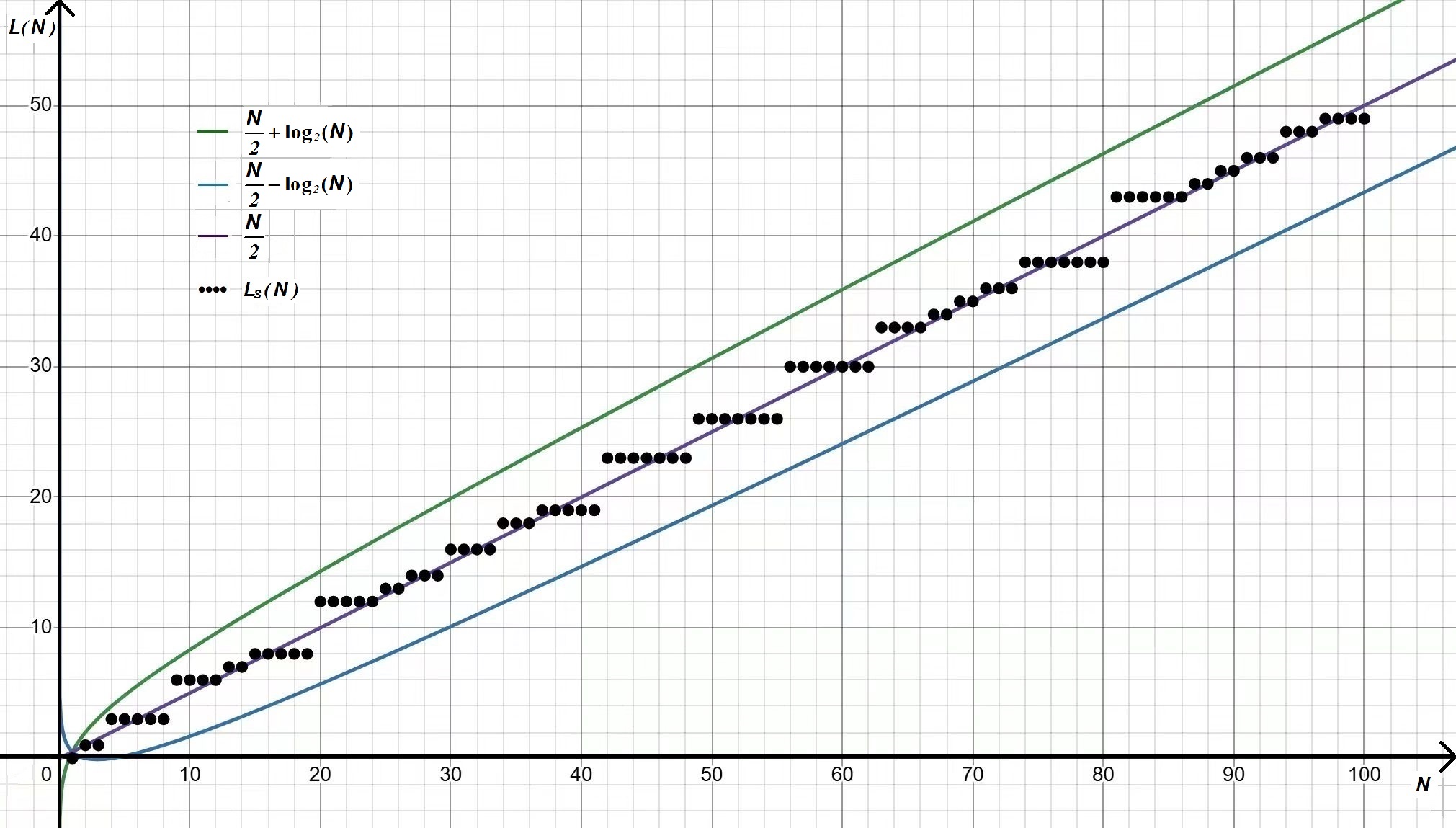} 
\end{center}
\caption{The $N$th linear complexity of $\mathcal{S}$ with $G_{\mathcal{S}}(2)=\sqrt{-7}$ for $N\le 100$}
\label{fig=-7}
\end{figure}

Similarly, our experimental data supports the conjecture that the $N$th linear complexity of both sequences $\mathcal{S}$ with $G_{\mathcal{S}}(2)=\sqrt{17}$ and $G_\mathcal{S}(2)=\sqrt{-7}$ is close to $\frac{N}{2}$, see Figs.\ref{fig=17} and \ref{fig=-7}.  

Undesirable features of automatic sequences are their small subword complexity, see \cite[Corollary 10.3.2]{AS2003} or \cite[Sect.\ 7]{MW2022}, 
and, under certain restrictions, their large $N$th well-distribution measure and $N$th correlation measure of order $k$, 
which are of worst possible order of magnitude $N$, which imply that they are statistically distinguishable from random sequences, see \cite{MS1998,MW2018CC,MW2022}.
In particular, for automatic sequences the subword complexity $p(n)$, that is the number of different subwords of length $n$ in the sequence, satisfies $p(n)=O(n)$, see for example \cite[Theorem 7.2]{MW2022}.
However, for irrational $2$-adic algebraic sequences we have 
$$\liminf_{n\rightarrow \infty} \frac{p(n)}{n}=\infty$$
by \cite[Theorem 1B]{AB2007} and these sequences behave better than automatic sequences with respect to the subword complexity, see also the survey \cite{B2015}.

Besides this result, we have neither detected any undesirable feature of $2$-adic algebraic sequences, yet, nor found any theoretical support that these sequences have desirable features with respect to the above-mentioned and related measures of pseudorandomness.

\subsection{A consequence of Roth's theorem}

Note that it follows from the analog of Roth's theorem for $2$-adic numbers that for any $\varepsilon>0$, any non-periodic sequence $\mathcal{S}$ satisfies for sufficiently large $N$,
\begin{equation}\label{roth}
\lambda_{\mathcal{S}}(N)> \frac{N}{2+\varepsilon},\quad N\ge N_0(\varepsilon).
\end{equation}
More precisely, for an irrational $\xi$ let $w_1(\xi)$ be the upper limit of the real numbers $w$ for which there exist infinitely many $(f,q)$ with 
$$\left|\xi-\frac{f}{q}\right|_2\le \max\{|f|,|q|\}^{-w-1}.$$
Then Roth's theorem, see \cite[Theorem~9.4]{B2004}, states $w_1(\xi)=1$.
Hence, for any $\varepsilon>0$ there are only finitely many $(f,q)$ with
$$\left|\xi-\frac{f}{q}\right|_2\le \max\{|f|,|q|\}^{-2-\varepsilon}$$
and for sufficiently large $N$ we have for all $(f,q)$,
\begin{equation}\label{heightbound}
2^{-N}\ge \left|\xi-\frac{f}{q}\right|_2>\max\{|f|,|q|\}^{-2-\varepsilon}.
\end{equation}
Applying this to $\xi=G_{\mathcal{S}}(2)$ we get 
$$2^{-N}>\Lambda_{\mathcal{S}}(N)^{-2-\varepsilon},\quad N\ge N_0(\varepsilon),$$
and \eqref{roth} follows.

Theorem \ref{Thm-bound} complements \eqref{roth} for small $N$.

\subsection{Real algebraic sequences}

The pseudorandomness of binary sequences $\mathcal{S}$ for which $2^{-1}G_{\mathcal{S}}(2^{-1})\in [0,1]$ is algebraic in $\mathbb{R}$ has also been studied, see for example \cite{SY2025}.
An analog of the $2$-adic complexity, the {\em rational complexity}, was introduced in \cite{VCJ2022}.
By Roth's theorem \cite[Theorem~2.1]{B2004} we get the analog of \eqref{heightbound} for the absolute value $|.|$. For any $(f,q)$ with $0<f<q$ and any non-rational $\xi$ we have 
$$2^{-N}\ge \left|\xi-\frac{f}{q}\right|>q^{-2-\varepsilon}$$
for sufficiently large $N$.

\section*{Acknowledgments} 
The authors would like to thank the anonymous referees and the associate editor for useful comments.

The work was written during a pleasant visit of Z.\ Chen to Linz. He wishes to thank the Chinese Scholarship Council for financial support and  RICAM of the Austrian Academy of Sciences for hospitality.

The work of Z.\ Chen was supported in part by the National Science Foundation of China under Grant No.\ 62372256 and by the Fujian Provincial Science Foundation of China under Grant No.\ 2023J01996. The research of A.\ Winterhof was funded in part by the Austrian Science Fund (FWF) [10.55776/PAT4719224].


\begin{thebibliography}{00}

\bibitem{AB2007}
B. Adamczewski, Y. Bugeaud.
On the complexity of algebraic numbers I. Expansions in integer bases. Ann. Math. 165(2): 547--565 (2007)

\bibitem{ACG2020}
B. Adamczewski, J. Cassaigne, M. Gonidec.
On the computational complexity of algebraic numbers: the Hartmanis-Stearns problem revisited. 
Trans. Amer. Math. Soc. 373(5): 3085--3115 (2020)
 

\bibitem{AHN2020}
J.-P. Allouche, G. Han, H. Niederreiter.
Perfect linear complexity profile and apwenian sequences. Finite Fields Their Appl. 68: 101761 (2020)

\bibitem{AS2003}
J.-P. Allouche, J. Shallit. Automatic Sequences: Theory, Applications, Generalizations. Cambridge
University Press, Cambridge, 2003.

\bibitem{BS1986}
Z. I. Borevich, I. R. Shafarevich. Number Theory.
 1st edition, vol.\ 20, Academic Press, New York (1986)


\bibitem{B2004} Y. Bugeaud. Approximation by Algebraic Numbers. Cambridge Tracts in Mathematics vol.~160,  Cambridge University Press, Cambridge, 2004.

\bibitem{B2015} Y. Bugeaud.
Expansions of algebraic numbers. Four faces of number theory, 31--75.
EMS Ser. Lect. Math.
European Mathematical Society (EMS), Z\"urich (2015)

\bibitem{BY2017}
Y. Bugeaud, J. Yao. Hankel determinants, Pad\'e approximations, and irrationality exponents for $p$-adic numbers. Annali di Matematica  196(3): 929--946 (2017) 

\bibitem{CCOW2024}
Z. Chen, Z. Chen, J. Obrovsky, A. Winterhof. Maximum-order complexity and $2$-adic complexity. IEEE Trans. Inf. Theory 70(8): 6060--6067 (2024)


\bibitem{CW2024}
Z. Chen, A. Winterhof. Probabilistic results on the $2$-adic complexity. Des. Codes Cryptogr.: to appear (2025). \url{https://doi.org/10.1007/s10623-025-01592-1} 

\bibitem{C1979}
G. Christol. Ensembles presque p\'{e}riodiques $k$-reconnaissables.  Theor. Comput. Sci. 9(1): 141--145 (1979)

 \bibitem{CKMR1980}
 G. Christol, T. Kamae, M. Mend\`{e}s France, G. Rauzy. Suites alg\'{e}briques, automates et substitutions.
Bull. Soc. Math. France 108(4): 401--419 (1980)

\bibitem{DPS2021}
J. Damien, P. Popoli, T. Stoll. Maximum order complexity of the sum of digits function in Zeckendorf base and polynomial subsequences.  Cryptogr. Commun. 13(5):   791--814 (2021)

\bibitem{D2012}
C. Diem. On the use of expansion series for stream ciphers. LMS J. Comput. Math. 15: 326--340 (2012) 

\bibitem{DMR2019}
M. Drmota, C. Mauduit, J. Rivat. Normality along squares. J. Eur. Math. Soc. 21(2): 507--548 (2019)

\bibitem{GM2020}
D. G\'omez-P\'erez, L. M\'erai.
Algebraic dependence in generating functions and expansion complexity. Adv. Math. Commun. 14(2): 307--318 (2020)

\bibitem{GMN2018}
D. G\'omez-P\'erez, L. M\'erai, H. Niederreiter.
On the expansion complexity of sequences over finite fields. IEEE Trans. Inf. Theory 64(6): 4228--4232 (2018)

\bibitem{GK2012}
M. Goresky, A. Klapper. Algebraic Shift Register Sequences. Cambridge University Press, Cambridge (2012)

\bibitem{G2024}
J. Gr\"{u}nberger. Linear Complexity and Continued Fractions of Families of Binary Sequences. Master Thesis, Johannes Kepler University Linz (2024)
\url{https://epub.jku.at/urn:nbn:at:at-ubl:1-72957}

\bibitem{HW2016}
R. Hofer, A. Winterhof.
Linear complexity and expansion complexity of some number theoretic sequences. International Workshop on the Arithmetic of Finite Fields -- WAIFI 2016, Lecture
Notes in Comput. Sci., vol.\ 10064: 67--74, Springer-Verlag, Cham (2016)

\bibitem{I2004}
R. S. Irving. Integers, Polynomials, and Rings: a Course in Algebra. Springer-Verlag, Berlin (2004)

\bibitem{K1984}
N. Koblitz. $p$-adic Numbers, $p$-adic Analysis and Zeta-Functions. 2nd edition, Springer-Verlag, Berlin (1984)

%\bibitem{MS1997} 
%C. Mauduit, A. S\'ark\"{o}zy. On finite pseudorandom binary sequences. I. Measure of pseudorandomness, the Legendre symbol. Acta Arith. 82(4): 365–377 (1997)

\bibitem{MS1998}
C. Mauduit, A. S\'ark\"{o}zy. On finite pseudorandom binary sequences II: The Champernowne, Rudin-Shapiro, and Thue-Morse sequences, a further construction. J. Number Theory 73(2): 256--276 (1998)
 
\bibitem{MNW2017}
L. M\'erai, H. Niederreiter, A. Winterhof.
Expansion complexity and linear complexity of sequences over finite fields. Cryptogr. Commun. 9(4): 501--509 (2017)
 
\bibitem{MW2022}
L. M\'{e}rai, A. Winterhof.
Pseudorandom sequences derived from automatic sequences. Cryptogr. Commun. 14(4): 783--815 (2022)

\bibitem{MW2018CC}
L. M\'{e}rai, A. Winterhof.
On the pseudorandomness of automatic sequences. Cryptogr. Commun. 10(6): 1013--1022 (2018)

\bibitem{MW2018JNT}
L. M\'{e}rai, A. Winterhof. On the $N$th linear complexity of automatic sequences. J. Number Theory 187: 415--429  (2018)



\bibitem{M2018}
C. M\"ullner. The Rudin-Shapiro sequence and similar sequences are normal along squares. Canad. J. Math. 70(5) : 1096--1129 (2018) 

\bibitem{MS2017}
C. M\"ullner, L. Spiegelhofer. Normality of the Thue-Morse sequence along Piatetski-Shapiro  sequences, II, Israel J. Math. 220 : 691--738 (2017)

\bibitem{N1990}
H. Niederreiter.
A combinatorial approach to probabilistic results on the linear complexity profile of random sequences. J. Cryptol. 2(2), 105--112 (1990)

\bibitem{N1999}
J. Neukirch.  
Algebraic Number Theory.
Grundlehren Math. Wiss., vol.\ 322.
Springer-Verlag, Berlin (1999)

\bibitem{O2023}
 J. Obrovsky. $2$-Adic Complexity and Related Measures of Pseudorandomness. Master Thesis, Johannes Kepler University Linz (2023)
\url{https://epub.jku.at/urn/urn:nbn:at:at-ubl:1-64657}

\bibitem{P2020}
P. Popoli. On the maximum order complexity of Thue-Morse and Rudin-Shapiro sequences along polynomial values. Unif. Distrib. Theory 15(2): 9-22 (2020)

\bibitem{SY2025}
A. Saito, A. Yamaguchi.
Accelerating true orbit pseudorandom number generation using Newton's method. 
\url{https://arxiv.org/abs/2502.16108}


\bibitem{S2015}
 L. Spiegelhofer.  Normality of the Thue-Morse sequence along Piatetski-Shapiro sequences. Q. J. Math. 66 : 1127--1138 (2015)

\bibitem{SW2019UDT}
Z. Sun, A. Winterhof. On the maximum order complexity of the Thue-Morse and Rudin-Shapiro sequence. Unif. Distrib. Theory 14(2): 33-42 (2019)

\bibitem{SW2019}
Z. Sun, A. Winterhof.
On the maximum order complexity of subsequences of the Thue-Morse and Rudin-Shapiro sequence along squares. Int. J. Comput. Math. Comput. Syst. Theory 4(1): 30-36 (2019)


\bibitem{VCJ2022}
M. Vielhaber, M. del Pilar Canales Chacón, S. Ceballos.
Rational complexity of binary sequences, F$\mathbb{Q}$SRs and pseudo-ultrametric continued fractions in $\mathbb{R}$.
Cryptogr. Commun. 14(2): 433--357 (2022).


\bibitem{W2023}
A. Winterhof.
Pseudorandom binary sequences: quality measures and number theoretic constructions. IEICE Trans. Fundam. Electron. Commun. Comput. Sci.
E106-A(12): 1452-1460 (2023)






\end{thebibliography}
\end{document}